\documentclass[10pt]{article}

\date{Universit\'e Paris-Sud}
\usepackage{amssymb,babel}
\usepackage{amsmath}
\usepackage{amsthm}
\usepackage{amsfonts}
\usepackage{amssymb,amsbsy,amsmath,amsfonts,amssymb,amscd}

\def \Prob {{\bf P}}
\newcommand{\R}{\mathbb{R}}

\newcommand{\Z}{\mathbb{Z}}

\def \H{\mathbb{H}}

\def \eps {{\varepsilon}}

\def \AA {{\mathcal A}}
\newtheorem {proposition}{Proposition}
\begin{document}

\title{
Conformal fields, restriction properties,\\
 degenerate
representations and SLE
\footnote {To appear in C.R. Acad. Sci. Paris, Ser. I Math.}
}

\author{
Roland Friedrich  \footnote {rolandf@ihes.fr}
\ \
\and   
Wendelin Werner
\footnote {wendelin.werner@math.u-psud.fr}%~$^{\text{b}}$
}
\maketitle
%%%%%%%%%%%%%%%%%%%%%%%%%%%%%%%%%%%%%%%%%%%%%%%%%%%%%%%%%%%%
%%%  Abstract  %%%
%%%%%%%%%%%%%%%%%%
%%%%%%%%%%%%%%%%%%%%%%%%%%%%%%%%%%%%%%%%%%%%%%%%%%%%%%%%%%%%
%%%  Abstract  %%%
%%%%%%%%%%%%%%%%%%
\begin{abstract}%
{We show how to relate the Schramm-Loewner Evolution
processes (SLE) to 
highest-weight representations of the Virasoro Algebra.
The restriction 
properties of SLE that have been recently derived in \cite {LSWr}
 play a crucial role. In this setup, various considerations
from conformal field theory can be interpreted and reformulated via SLE. 
This enables to make a concrete link between the two-dimensional 
discrete critical systems from statistical physics and conformal 
field theory. 
}
\end{abstract}
%%%%%%%%%%%%%%%%%%%%%%%%%%%%%%%%%%%%%%%%%%%%%%%%%%%%%%%%%%%%
%%%  Main text (in English)  %%%
\section {Introduction}
Conformal field theory
 has been remarkably successful in predicting the 
critical behaviour of two-dimensional 
systems from statistical physics (see  \cite {BPZ,Ca1} and
for instance the compilation of papers in \cite {ISZ}).
One of the starting points \cite {Po,BPZ} for this theory 
is that each system should be (in its scaling limit) correspond
to a conformal field (these fields are classified according to their central
charge).
The behaviour of the system should then be described
by critical exponents which physicists identify as highest-weights of 
certain degenerate representations of infinite-dimensional Lie Algebras
(this motivated also many mathematical papers on representation
theory, one should 
probably cite here at least all 
the papers reprinted in \cite {GO}, for the 
representation theory of the Virasoro Algebra and 
background, see \cite {KR,Ka}).
This applies for instance to Ising and Potts models, percolation, 
self-avoiding walks. 
However, as for instance pointed out in \cite {LPS}, many fundamental 
features  have remained unclear at least for mathematicians.
For instance, the actual relation between the discrete system and the 
fields ie. the meaning of the field in terms of the discrete 
system, the interpretation of the equations  
leading to the highest-weight representation were rather mysterious. 

In 1999, Oded Schramm \cite {S1}
defined a one-parameter family of random curves based
on Loewner's differential equation, SLE$_\kappa$
indexed by the positive real parameter $\kappa$ (SLE stands
for Stochastic -or Schramm- Loewner Evolution).
These random curves are the only ones which
combine conformal invariance and a Markovian-type property. Provided
that the scaling limit of interfaces in the above-mentioned models
exist and are conformally invariant, then  the limiting objects
must therefore  be one of the
SLE$_\kappa$ curves.
This has now been rigourously proved in some cases
(critical site
percolation on the triangular lattice
\cite {Sm} and uniform spanning trees 
\cite {LSWlesl}).
For a general discussion of the conjectured relation between the
discrete models and SLE, see \cite {RS}.
In this SLE setting, the 
critical exponents are simply principal eigenvalues of some differential operators
\cite {LSW1,LSW2,LSW3,LSW5}.
This led to a complete mathematical proof for the value of
critical exponents for those models that have been proved to 
be conformally invariant. 
In order to confirm  rigorously the conjectures for the other models,
the missing step is to derive their
conformal invariance.

It is therefore natural to think that SLE should be related to 
conformal field theory (as for instance hinted at in \cite {BB}) and 
to 
highest-weight representations of the Virasoro Algebra.
Our goal in the present note is to explain that it is indeed the
case. 
In order to make this exposition as clear 
and rigorous as possible (and also to keep this note short),
we will restrict ourselves to the simplest case: The ``boundary
behaviour'' of SLE$_{8/3}$ that corresponds in the field theory 
language to a zero central charge and conjecturally corresponds to the 
scaling limit of the half-plane self-avoiding walk \cite {LSWSAW}.
The other cases (non-zero central charge, behaviour in the bulk)
will be detailed in forthcoming papers.
 
\section {SLE facts}

The chordal SLE$_\kappa$ curve $\gamma$
is characterized as follows: 
The conformal maps
$g_t$ from $\H \setminus \gamma [0,t]$ onto $\H$ 
such that $g_t (z) = z + o(1) $ when $z \to \infty$
solve 
the ordinary differential equation 
$\partial_t g_t (s) =
2 /( g_t(z) - W_t)$ 
(and is started from $g_0(z) = z$),
where
$W_t = \sqrt {\kappa} \beta_t$ (here and in the 
sequel, $(\beta_t, t \ge 0)$ is a standard real-valued 
Brownian motion with $\beta_0=0$).
In other words, $\gamma_t$ is precisely the
point such that $g_t(\gamma_t) = W_t$.
See e.g. 
\cite {LSW1, RS} for the definition and 
properties of  SLE, or \cite {Lin, Wstf} for reviews.

In the recent paper \cite {LSWr}, it is shown how to
measure the distortion of the law of an SLE when
its image is mapped conformally from a subdomain onto some other
domain.
In particular, for the special value $\kappa =8/3$, the SLE curve $\gamma$
(which is conjectured to be the scaling limit of a certain measure on
self-avoiding curves, see \cite {LSWSAW} for a discussion of this conjecture
and some of its consequences -- this conjecture is conforted by
numerical simulations \cite {Kenn}) has the
{\em conformal restriction
property} that we now briefly describe:

Suppose that $H$ is a
simply connected open subset of the upper-half plane
$\H$
such that $\H \setminus H$ is bounded and bounded away from $0$.
Assume also that $\gamma$ is chordal SLE$_{8/3}$.
Then, the law of $\gamma$ conditioned
to remain in $H$ is identical to the law of $\Phi(\gamma)$ where $\Phi$
is the conformal map from $\H$ onto $H$ that fixes the two boundary 
points $0$ and $\infty$ and $\Phi (z) \sim z$ as $z \to \infty$.
Actually \cite {LSWr}, SLE$_{8/3}$
is the unique simple random curve with this
property. As argued towards
the begining of the paper  \cite {LSWr},
 if $\gamma$ satisfies the restriction
property, then there exists $\alpha >0$ such that
\begin {equation}
\label {rest}
 P [ \gamma \subset H ] =   \Phi'(0)^{-\alpha}.
\end {equation}
It is proved in \cite {LSWr}
 that the value of $\alpha$ corresponding to SLE$_{8/3}$
is $5/8$ (but we will not use this 
fact in this note, as one of our goals is to explain why 
one can recover it from algebraic considerations as was 
for instance predicted in \cite {DS}).

\section {Boundary behaviour}
When $\kappa <8$ and 
 $\eps \to 0$,
the probability that chordal SLE$_\kappa$
intersects the $\eps$-neighbourhood of the real point $x \not=0 $
 can be shown 
to decay (up to a multiplicative constant) like
$\eps^{s}$ when $\eps \to 0$,
where 
$s = (8/ \kappa) - 1$.
More generally, one defines 
for $x_1, \ldots, x_n \in \R \setminus \{ 0 \}$,
$$
B_n (x_1, \ldots , x_n)
=\lim_{\eps_1 , \ldots, \eps_n \to 0}
\eps_1^{-s} \ldots \eps_n^{-s}
 \Prob [ \cap_{j=1}^n
\{ \gamma \cap [x_j, x_j + i \eps_j \sqrt {2}]  \not= \emptyset  \} ]
$$
(the choice of these vertical slits enables 
one to have simple renormalizing constants later on).
It is easy to see using (\ref {rest}) that if the restriction 
property holds, then $s=2$ and $B_1 (x) = \alpha / x^2$.
The functions $B_n$ correspond to correlation functions in conformal
field theory. They have the following properties:

\begin {proposition}
\null
\begin {enumerate}
\item {Ward-type identity:}
If the restriction property holds (with the exponent $\alpha$),
 then for all $n \ge 1$,
\begin {eqnarray}
\nonumber
\lefteqn{ B_{n+1} ( x, x_1, x_2, \ldots, x_n )
\ = \  \frac {\alpha}{x^2} B_n (x_1, \ldots , x_n)} \\
&&
%B_{n+1} ( x, x_1, \ldots, x_n ) 
%=
%\left[ 
%\frac {\alpha}{x^2}
-\sum_{j=1}^n \left\{ (\frac {1}{x_j-x} + \frac 1x ) \partial_{x_j}
  - \frac {2}{(x_j-x)^2} \right\} %\right]
 B_n (x_1, \ldots , x_n).
\label {e.ward}
\end {eqnarray}
\item {Evolution equation:}
\begin {equation}
\label {e.bdy}
-2 s \left( \sum_{j=1}^n \frac {1}{x_j^2} \right) B_n
+ \left( \sum_{j=1}^n \frac {2}{x_j} \partial_{x_j} \right)
B_n
+ \frac {\kappa}{2}  \left( \sum_{j=1}^n \partial_{x_j} \right)^2
B_n = 0,
\end {equation}
\end {enumerate}
\end {proposition}

\noindent
{\bf Proof (sketch).}
The first identity is in fact a direct consequence of the restriction 
property:
Suppose that the real numbers $x_1, \ldots, x_n$ are fixed and let
us focus on the event  that $\gamma$ does intersect all the
slits  $[x_j, x_j + i \eps \sqrt {2}]$.  Let us also choose another
point $x \in \R$ and a small $\delta$.  Now, either $\gamma$ avoids $[x ,
x+ i \delta \sqrt {2}]$ or it does also hit it. 
The probability of the first event can be expressed using the
restriction property, while the second one involves the 
$(n+1)$-points functions $B_{n+1}$. The sum of the two probabilities 
remains independent of $\delta$. Looking at the $\delta^2$
term in its expansion yields (\ref {e.ward}). 
Let us stress that these identities together with the knowledge of
$B_1 (x_1) = \alpha  /x_1^{2}$ fully determine all the functions
$B_n$. One can also write $B_0=1$ as a function of $0$ variables.

The Markovian property of the SLE curves (which is a straightforward
consequence of the stationarity of the increments of the Brownian motion $
\beta$)
goes as follows:
The law of $\gamma [t, \infty)$ given $\gamma [0,t]$ is
the image of an independent copy $\tilde \gamma$ of $\gamma$ under a conformal map
from $\H$ on $\H \setminus \gamma [0,t]$ which maps $0$ onto $\gamma(t)$ and
$\infty$ onto itself.
In other words,
the conditional law
 of $g_t ( \gamma [t, \infty)) - W_t $ is independent of $\gamma [0,t]$
and identical to that of the initial SLE.
This shows immediately that for all fixed reals $x_1, \ldots, x_n$,
 the processes
$$ |g_t'(x_1)|^{s} \ldots |g_t'(x_n)|^{s}
B_n ( g_t (x_1) - W_t , \ldots , g_t (x_n) - W_t )$$
are local martingales.
It\^o's formula then yields immediately (\ref {e.bdy}).
\qed

\section {Representations}

It is interesting to focus on the asymptotic expansion of
$B_{n+1} (x, x_1, \ldots, x_n)$ when $x \to 0$.
It is natural to define the operators
$$
{\cal L}_{-N}
=
\sum_{j}
\left\{ - x_j^{1-N} \partial_j + 2(N-1)x_j^{-N} \right\}
.$$
Note that these operators satisfy the commutation relation
\begin {equation}
\label {commute}
[{\cal L}_N, {\cal L}_{M}]
= (N-M) {\cal L}_{N+M}.
\end {equation}
The Algebra generated by vectors satisfying this commutation 
relation will be denoted by $\AA$ 
and is often 
viewed as the Lie algebra of vector fields 
of the circle i.e. 
the algebra generated by $e_N = -  z^{N+1} d/dz, N \in \Z$ (which 
satisfy the same commutation relation).
 One can rewrite the previous Proposition using the 
operators ${\cal L}$:
The Ward-type identity becomes
\begin {equation}
\label {ward2}
B_{n+1} (x, x_1, \ldots, x_n )
=
\frac {\alpha}{x^2}   B_n (x_1, \ldots, x_n)
+ \sum_{N \ge 1} x^{N-2} {\cal L}_{-N} B_n  (x_1, \ldots, x_n )
\end {equation}
if $|x| < \min (x_1, \ldots, x_n)$.
The evolution equation becomes (if $s=2$ i.e. if the restriction
property holds) 
\begin {equation}
\label {L2}
(\frac {\kappa}{2}
{\cal L}_{-1}^2
-2 {\cal L}_{-2}) B_n
= 0
.\end {equation}

It is easy to make a little computation to see  that 
$\kappa=8/3$ and $\alpha = 5/8$ if
the previous Proposition holds for some family of functions $B_n$.
This computation has some similarities with the computation that 
shows that a highest-weight representation of $\AA$ 
that is degenerate at level two has a highest weight equal to $5/8$.
We now show that this is not just a similarity.
Note that scaling implies that ${\cal L}_0 B_n = 0$, 
so that the representation
is not simply given in terms of the ${\cal L}$'s.

Here is one way to construct a 
highest-weight representation of $\AA$.
Define the vector
${B}  = ( B_0, B_1, \ldots )$.
Suppose that $w = (w_0, w_1, \ldots)$ is such that
$w_n$ is a (rational) function of the $n$ variables $x_1, \ldots , x_n$.
We define the operators $l_N$ is such a way that
\begin {equation}
\label {Laurent}
w_{n+1} (x, x_1, \ldots, x_n)
= \sum_{N \in \Z} x^{N-2} (l_{-N}(w))_n (x_1, \ldots, x_n )
\end {equation}
when $x \to \infty$.
In other terms the function of $n$ variables in $l_{-N} (w)$ is the
$x^{N-2}$ term in the Laurent  expansion of
$w_{n+1} (x, x_1, \ldots, x_n)$.
Equation (\ref {ward2})
shows that
\begin {equation}
\label {hw}
l_N (B) = \left\{
 \begin {array}{l@{\hbox { if }}l}
 (0, 0, \ldots ) & N>0 \\
(\alpha B_0, \alpha B_1, \ldots )  & N=0 \\
({\cal L}_N B_0, {\cal L}_N B_1, \ldots)  & N<0
\end {array}
\right.
\end {equation}
We then define the vector space $V$ generated by $B$ and all the vectors
$ l_{N_1} \cdots l_{N_r} (B)$
for $N_1 \le  \cdots \le N_r < 0$.
It is then possible to prove (without using the evolution equation)
that:
\begin {proposition}
\label {algebra}
\begin {itemize}
\item
For $N_1 \le  \cdots \le N_r < 0$,
the function of $n$ 
variables in $l_{N_1} \cdots l_{N_r} (B)$ is 
${\cal L}_{N_1} \cdots {\cal L}_{N_r} B_n$.
\item
$V$ is stable under all $l_N$'s (for all $N \in \Z$, not only for $N \le 0$).
\item
$[l_M  , l_N] w = (M-N) l_{M+N} w $
for all $w \in V$ and all $N,M \in \Z$ (note again that $M$ and $N$ can be 
positive).
\end {itemize}
\end {proposition}

\medbreak
\noindent 
{\bf Proof (sketch).}
The first statement can be proved inductively using (\ref {ward2}).
It implies immediately the third statement for all negative $N,M$.
In order to obtain the commutation relation 
for all $N,M$, and due to the fact that $l_N (B) =0 $ for all positive $N$,
it in fact suffices to check it if only $M$ is positive. 
This can then be done again thanks to (\ref {ward2}) and the 
first statement.
Finally, note that the second statement is a consequence of the third one
and of the fact that $l_N (B) = 0$ if $N >0$.
\qed

\medbreak

Equation (\ref {hw}) shows that this representation of the 
algebra $\AA$ 
generated by the $l_N$'s is a 
highest-weight representation with highest weight $\alpha$.
Equation (\ref {L2}) means that this representation is degenerate at level 2,
and it is easy to verify that this implies that the highest weight $\alpha$
is equal to $5/8$ (see \cite {GO}).
In other words, the fact that one obtains a representation of  
$\AA$ is a consequence of the restriction property, and 
its degeneracy follows from the Markovian property.

\section {Generalizations}

Analogous arguments can be used
 to relate other chordal SLE's to  
degenerate 
highest-weight representations of the Virasoro algebra (which is the
central extension of $\AA$).
When $\kappa \not= 8/3$, SLE
does not satisfy the conformal restriction
property. However, it is possible to quantify
how much the property fails to be true
and this \cite {LSWr} involves the integral
of a certain constant $c(\kappa)$ times
the Schwarzian derivative of some
conformal maps (see also \cite {LSWr,LW} for 
another interpretation).
This constant $c(\kappa)$ will 
turn out to be exactly the central charge of the 
corresponding representation.
The highest-weight is the exponent of the 
restriction measure that is naturally associated to the SLE
via the correspondance derrived in \cite {LSWr}.
This will be detailed in a forthcoming paper.

One can also make a similar analysis for the SLE behaviour in the bulk, and
make the link to Ward-type
 identities and  representation theory. However, in this 
case, there are some additional questions to solve if 
one looks for a completely rigorous link.
For instance,  the
very  
definition of the ``correlation function''
(see \cite {Be3} to see the difficulty of 
estimating precisely the probability that the SLE paths passes 
in the neighbourhood of $n$ 
given points in the upper half-plane) is problematic.

%%%%%%%%%%%%%%%%%%%%%%%%%%%%%%%%%%%%%%%%%%%%%%%%%%%%%%%%%%%
%%%  Acknowledgements  %%%
%%%%%%%%%%%%%%%%%%%%%%%%%%

\medbreak
\noindent
{\bf Acknowledgements.}
This work ows a lot to stimulating discussions with 
Vincent Beffara and Yves Le Jan. 
Thanks are also due to Greg Lawler and Oded Schramm,
in particular because of the instrumental role played by the ideas 
developed in \cite {LSWr}.
R.F. acknowledges support and hospitality of IHES.
%%%%%%%%%%%%%%%%%%%%%%%%%%%%%%%%%%%%%%%%%%%%%%%%%%%%%%%%%%%%
%%%  Bibliography  %%%
%%%%%%%%%%%%%%%%%%%%%%%

------------------

Laboratoire de Math\'ematiques

Universit\'e Paris-Sud

91405 Orsay cedex, Fracne

 \end {document}